\documentclass[11pt]{amsart}

\usepackage{amssymb,latexsym,vmargin,graphicx,color}
\usepackage[colorlinks,urlcolor=blue,unicode,bookmarksopen,pdfstartview=FitH]{hyperref}

\newtheorem{theorem}{Theorem}

\newtheorem{lemma}{Lemma}

\theoremstyle{definition}

\newcommand{\beql}[1]{\begin{equation}\label{#1}}
\newcommand{\eeq}{\end{equation}}
\newcommand{\comment}[1]{}

\newcommand{\Abs}[1]{{\left|{#1}\right|}}

\newcommand{\Prob}[1]{{{\bf{Pr}}\left[{#1}\right]}}
\newcommand{\Set}[1]{{\left\{{#1}\right\}}}

\newcommand{\RR}{{\mathbb R}}
\newcommand{\CC}{{\mathbb C}}
\newcommand{\ZZ}{{\mathbb Z}}

\newcommand{\ft}[1]{\widehat{#1}}

\newcounter{rem}
\begin{document}

\title[The discrepancy of a needle on a checkerboard]{The discrepancy of a needle on a checkerboard}

\author[M. Kolountzakis]{Mihail N. Kolountzakis}
\address{Department of Mathematics, University of Crete, Knossos Ave., GR-714 09, Iraklio, Greece}
\email{kolount@gmail.com}

\date{October 2007}

\thanks{
Supported by the Greek research program ``Pythagoras 2'' (75\% European funds
and 25\% National funds) and by INTAS 03-51-5070 (2004)
({\em Analytical and Combinatorial Methods in Number Theory and Geometry}).
}

\begin{abstract}
Consider the plane as a checkerboard, with each unit square
colored black or white in an arbitrary manner.
We show that for any such coloring there are straight line segments, of arbitrarily large length,
such that the difference of their white length minus their black length,
in absolute value, is at least the square root of their length, up to a multiplicative constant.
For the corresponding ``finite'' problem ($N \times N$ checkerboard) we also prove that we can
color it in such a way that the above quantity is at most $C \sqrt{N \log N}$,
for any placement of the line segment.
\end{abstract}

\maketitle

\section{Introduction}

In this paper we answer a question posed to us by
\href{http://users.uoa.gr/~ppapazog/}{P. Papasoglu} \cite{papasoglu}:
\begin{quote}
Consider the plane as a checkerboard, with each unit square
$$
[m, m+1) \times [n,n+1),\ \  m,n \in \ZZ,
$$
colored black or white.
Is it possible that there is such a coloring and a finite constant $M$ such that
for any line segment $I$ placed on the checkerboard the difference of its white length minus
its black length is, in absolute value, at most $M$?
\end{quote}
We show below that the answer is negative.
\begin{theorem}\label{th:infinite}
Suppose each unit cell $[i,i+1) \times [j,j+1),\ i,j \in \ZZ,$ in the plane is painted black or white.
Then there exist arbitrarily long line segments $I$ such that the difference of the black part of $I$ minus
its white part is at least $C \sqrt{\Abs{I}}$ in absolute value, where $C>0$ is an absolute constant.
\end{theorem}
Our approach is Fourier analytic and Theorem \ref{th:infinite}
follows in an obvious way from the ``finite'' Theorem \ref{th:line-discrepancy} below, by setting $z_{ij}=1$ if the cell
$[i,i+1) \times [j,j+1)$ is black and $z_{ij}=-1$ if the cell is white. 
\begin{theorem}\label{th:line-discrepancy}
Suppose the function $f:\RR^2 \to \CC$ is zero outside
the square $[0,N]\times[0,N]$ and is constant and equal to $z_{ij}$ on each of
the squares $(i,j)+[0,1)^2$, $0\le i, j<N$.
Then there is a straight line $S$ such that
$$
\Abs{\int_S f} \ge C N^{-1/2} \left(\sum_{i,j=0}^{N-1} \Abs{z_{ij}}^2\right)^{1/2}, 
$$
where $C$ is a positive constant.
\end{theorem}

It is possible to give an almost matching upper bound for Theorem \ref{th:line-discrepancy}.
In Theorem \ref{th:upper-bound} we show how to color a $N\times N$ checkerboard so that all line segments
placed on the checkerboard have white-over-black excess (discrepancy) bounded by $C \sqrt{N \log N}$.
\begin{theorem}\label{th:upper-bound}
For each $N$ there exists a function $f:\RR^2\to\Set{0,\pm1}$ which is
zero outside the square $[0,N]\times[0,N]$ and is constant
and equal to $1$ or $-1$ on each of the squares $(i,j)+[0,1)^2$, $0\le i, j<N$,
and is such that on any straight line segment $I$ we have
\beql{segmentbound}
\Abs{\int_I f} \le C \sqrt{N \log N}.
\eeq
\end{theorem}
Unfortunately this does not translate to a coloring of the infinite checkerboard so that the discrepancy of
any line segment $I$ is $o(\Abs{I})$. This problem we leave open.

The question dealt with in this paper falls naturally into the subject of geometric discrepancy \cite{matousek}.
In this research area there is usually an underlying measure $\mu$ as well as a family ${\mathcal F}$ of
allowed subsets of Euclidean space, on which the measure $\mu$ is evaluated and upper and lower bounds
are sought on the range of $\mu$ on ${\mathcal F}$.
The most classical case is that where $\mu$ is a normalized collection of points masses in the unit square minus Lebesgue
measure and ${\mathcal F}$ consists of all axis-aligned rectangles in the unit square.
Usually the underlying measure $\mu$ has an atomic part (point masses) and the family ${\mathcal F}$
consists of ``fat'' sets.
In the problem we are studying here the measure $\mu$ has no atomic part (it is absolutely continuous)
and the collection ${\mathcal F}$ consists of all straight line segments, which may be considered thin sets,
and, strictly speaking, $\mu$ is $0$ on these sets.
We believe that the closest work in the bibliography to this paper is that of Rogers \cite{rogers} where
the measure $\mu$ is the same as here but the family ${\mathcal F}$ consists not of straight line segments
but of thin strips.
There does not seem to exist a connection of our results with those of \cite{rogers}.

\section{Proofs}

The definition of the Fourier Transform that we use is
$$
\ft{f}(\xi) = \int_{\RR^d} f(x) e^{-2\pi i \xi\cdot x}\,dx
$$
for $f \in L^1(\RR^d)$ and $\xi \in \RR^d$.

We use the letter $C$ in this paper as an absolute positive constant, not the same in all its occurrences.

\begin{proof}(of Theorem \ref{th:line-discrepancy})\\
For a straight line $L$ through the origin let us denote by $\pi_L f$ (the projection of $f$ onto $L$)
the function of $t \in \RR$ given by
$$
\pi_Lf (t) = \int_{\RR} f(tu+su^{\perp}) \,ds,
$$
where $u$ is a unit vector along $L$ and $u^\perp$ is a unit vector orthogonal to $u$.
It is well known and easy to see that the Fourier Transform of $\pi_Lf$ is equal to the restriction
of the (two-dimensional) Fourier Transform of $f$, $\ft f$, on $L$:
$$
\ft{\pi_Lf}(\xi) = \ft{f}(\xi u),\ \ \xi\in\RR.
$$
Write $M = \sup_{L,t} \Abs{\pi_Lf(t)}$ where the supremum is taken over all lines $L$ through the origin
and real numbers $t$.
We aim to prove that
$$
M>C N^{-1/2} \left(\sum_{i,j=0}^{N-1} \Abs{z_{ij}}^2\right)^{1/2},
$$
for some constant $C>0$.

Since the support of $f$ has diameter $\le C N$ it follows from Parseval's equality that
\beql{l2-line}
\int_{\RR} \Abs{\ft{f}(tu)}^2\,dt = \int_{\RR} \Abs{\pi_Lf(t)}^2\,dt \le CM^2 N.
\eeq
It also follows from Parseval's equality that
\beql{l2-plane}
\sum_{i,j=0}^{N-1} \Abs{z_{ij}}^2 = \int_{\RR^2} \Abs{f}^2 = \int_{\RR^2}\Abs{\ft{f}}^2.
\eeq
\begin{lemma}\label{lm:decay}
For sufficiently large $A>0$ we have
$$
\int_{\Abs{\xi}>A} \Abs{\ft{f}(\xi)}^2\,d\xi \le \frac{C}{A} \sum_{i,j=0}^{N-1} \Abs{z_{ij}}^2.
$$
\end{lemma}
\begin{proof}
If $Q=[0,1]^2$ then by a simple calculation
$$
\ft{\chi_Q}(\xi_1,\xi_2) = e^{2\pi i(\xi_1+\xi_2)/2}\frac{\sin\pi\xi_1}{\pi\xi_1} \frac{\sin\pi\xi_2}{\pi\xi_2}
$$
and, therefore, with $\xi=(\xi_1, \xi_2)$ we have
$$
\Abs{\ft{\chi_Q}(\xi)}^2 \le \frac{C}{(1+\Abs{\xi_1})^2 (1+\Abs{\xi_2})^2}.
$$
If $G = \Set{(m, n): 0\le m,n < N}$ then
$$
f = \chi_Q * \sum_{p \in G} z_p \delta_p,
$$
where $\delta_p$ denotes the unit point mass at $p$ and $z_p$ is the value of $f$
in the cell $p+[0,1)^2$.
It follows that $\ft{f}(\xi) = \ft{\chi_Q}(\xi) \phi(\xi)$, where
the trigonometric polynomial $\phi(\xi)$ is given by
$$
\phi(\xi) = \sum_{p \in G} z_p e^{2\pi i p\cdot \xi}.
$$
The function $\phi$ is $\ZZ^2$-periodic and $\int_Q \Abs{\phi}^2 = \sum_{i,j=0}^{N-1} \Abs{z_{ij}}^2$.
We have
\begin{eqnarray*}
\int_{\Abs{\xi}>A} \Abs{\ft{f}(\xi)}^2\,d\xi &\le&
 \int_{\Abs{\xi_1}>0.1 A\mbox{\ or\ }\Abs{\xi_2}>0.1 A} \Abs{\ft{f}(\xi)}^2\,d\xi \\
 &\le& \int_{\Abs{\xi_1}>0.1 A \mbox{\ or\ } \Abs{\xi_2}>0.1 A} \frac{C\Abs{\phi(\xi_1,\xi_2)}^2}{(1+\Abs{\xi_1})^2 (1+\Abs{\xi_2})^2} \,d\xi \\
 &\le& \sum_{m>0.1 A,\mbox{\ or\ }n>0.1 A} \frac{C}{(1+\Abs{m})^2(1+\Abs{n})^2} \int_{Q+(m,n)}\Abs{\phi}^2\\
 &=& \int_Q\Abs{\phi}^2 \cdot \sum_{m>0.1 A,\mbox{\ or\ }n>0.1 A} \frac{C}{(1+\Abs{m})^2(1+\Abs{n})^2}\\
 &\le& \sum_{i,j=0}^{N-1} \Abs{z_{ij}}^2 \cdot \frac{C}{A},
\end{eqnarray*}
which proves the lemma.
\end{proof}

It follows from Lemma \ref{lm:decay} and \eqref{l2-plane} that, for sufficiently large $A>0$, we have
$$
\frac{1}{2} \sum_{i,j=0}^{N-1} \Abs{z_{ij}}^2 < \int_{\Abs{\xi}<A} \Abs{\ft{f}(\xi)}^2 
$$
Rewriting the right hand side in polar coordinates and using \eqref{l2-line} we get
\begin{eqnarray*}
\frac{1}{2} \sum_{i,j=0}^{N-1} \Abs{z_{ij}}^2 &\le& A \int_{u \in S^1} \int_{\Abs{t}<A} \Abs{\ft{f}(tu)}^2 \\
 &\le& C A M^2 N, 
\end{eqnarray*}
which proves the desired inequality.
\end{proof}

\begin{proof} (of Theorem \ref{th:upper-bound})\\
The proof is probabilistic. We choose the value of $f$ in each cell $(i,j)$,
$0\le i, j <N$, to be $\pm1$ with equal probability and independently
of other cells.

Fix a line segment $I$ which is contained in $[0,N] \times [0,N]$ and let $\ell_1,\ldots,\ell_k$ be the lengths
of the intersections of $I$ with each of the cells it intersects.
Then $\int_I f = \sum_{j=1}^k \epsilon_j \ell_j$, where $\epsilon=\pm 1$
independently and uniformly.

We use the standard estimate (see for instance \cite[Appendix A]{alon-spencer}):
$$
\Prob{\Abs{\int_I f} > \lambda \sigma} \le C e^{-C\lambda^2}
$$
where $\sigma^2 = \sum_{j=1}^k \ell_j^2$.
Since the $\ell_j$ are bounded by $\sqrt2$ it follows that
$$
\sigma \le C \sqrt{\sum_{j=1}^k \ell_j} \le C \sqrt\Abs{I} \le C \sqrt N,
$$
hence, using $\lambda=C\sqrt{\log N}$ we obtain
\beql{estimate}
\Prob{\Abs{\int_I f} \ge C \sqrt{N \log N}} \le C N^{-K},
\eeq
where we can choose the positive number $K$ as large as we please, by choosing
the appropriate constant $C$ in the left hand side of \eqref{estimate}.

Our next step is to apply the bound \eqref{estimate} to a well chosen set ${\mathcal S}$
of line segments such that the validity of \eqref{segmentbound} for the elements of ${\mathcal S}$
implies the validity of
\eqref{segmentbound} for all segments (but with a larger constant). It is important that
$\Abs{\mathcal S}$ must not grow faster than a fixed power of $N$, which will allow us
to use \eqref{estimate} simultaneously for all segments in ${\mathcal S}$.

A straight line segment $I$ is determined by its endpoints, call them $a_I$ and $b_I$.
Our set ${\mathcal S}$ will contain all segments that are determined by
a grid of points of spacing $N^{-10}\times N^{-10}$.
We will show that for any line segment $I$ there is a line segment $J \in {\mathcal S}$
such that  
\beql{near-segments}
\Abs{\int_I f - \int_J f} \le 1.
\eeq
First, we assume that the line segment $I$ is not nearly horizontal or nearly vertical.
That is we assume that the line segment $I$ forms an angle between $N^{-1}$ and $\frac{\pi}{2}-N^{-1}$
with the $x$-axis.
In this case it is obvious that if the endpoints of $J$ are $O(N^{-10})$ apart from those of $I$
then, for any cell $K$ that $I$ intersects at length $\ell$, the intersection of $J$
with $K$ has length $\ell+O(N^{-5})$, from which \eqref{near-segments} clearly follows.

Let us now deal with the case of $I$ being nearly horizontal.
(The case of nearly vertical segments is treated similarly.)
In this case the segment $I$ is either entirely contained in one horizontal
strip
\beql{strip}
[0,N]\times[i,i+1]
\eeq
or is contained in two successive such strips.
The latter case can be treated by breaking up $I$ into two segments each of which is contained
in one strip, so we assume that $I$ is contained in strip \eqref{strip}.
Take then $J$ to be a member of ${\mathcal S}$ whose endpoints are $O(N^{-10})$ apart from those
of $I$ and which is also contained in strip \eqref{strip}, and \eqref{near-segments} follows.

Finally, by an application of \eqref{estimate} with $K=100$, we conclude, that with
positive probability, for all $I \in {\mathcal S}$ we have $\Abs{\int_I f} \le C \sqrt{N \log N}$.
Hence, the latter inequality is true with high probability for all line segments $I$.
\end{proof}


\end{document}